\documentclass{SCIYA2017enOL}
\online
\begin{document}

\ensubject{fdsfd}%二级学科

%%%%%%%%%%%%%%%%%%%%%%%%%%%%%%%%%%%%%%%%%%%%%%%%%%%%%%%
%%% Authors do not modify the information below
%%% 作者不需要修改此处信息
%%% 有专题名称时, 将第一行的{}注释掉, 使用第二行
\ArticleType{ARTICLES}%栏目
%\SpecialTopic{Progress of Projects Supported by NSFC}%专题
%\SubTitle{Dedicated to Professor Yang Lo on the Occasion of his {\rm 70}th Birthday}%专刊说明
\Year{2017}
\Month{January}%
%\Vol{60}
%\No{1}
\BeginPage{1} %
%\DOI{10.1007/s11425-000-0000-0}
\ReceiveDate{January 1, 2017}
\AcceptDate{January 1, 2017}
%\OnlineDate{January 1, 2017}
%%%%%%%%%%%%%%%%%%%%%%%%%%%%%%%%%%%%%%%%%%%%%%%%%%%%%%%

%%% title: 标题
%%%   \title{title}{title for citation}
\title[]{Pseudo asymptotically periodic solutions for fractional integro-differential neutral equations\thanks{It has been accepted for publication in SCIENCE CHINA Mathematics.}}
{Pseudo asymptotically periodic solutions for fractional integro-differential neutral equations}%%后边花括号是眉题

%%% Corresponding author: 通信作者
%%%   \author[number]{Full name}{{email@xxx.com}}
%%% General author: 一般作者
%%%   \author[number]{Full name}{}
\author[1,2]{YANG Min}{yangm58@mail2.sysu.edu.cn}
\author[1,$\ast$]{WANG Qiru}{mcswqr@mail.sysu.edu.cn}

%%% Author information for page head. 页眉中的作者信息
%%% 若此处指定以此处为准, 否则直接调用author信息
\AuthorMark{Yang M}

%%% Authors for citation. 首页引用中的作者信息
%%% 若此处指定以此处为准, 否则直接调用author信息
\AuthorCitation{Yang M, Wang Q R}

%%% Address. 地址
%%%   \address[number]{Address, City {\rm Postcode}, Country}
\address[1]{School of Mathematics, Sun Yat-sen University, Guangzhou {\rm 510275}, China}
\address[2]{School of Mathematics, Taiyuan University of Technology , Taiyuan {\rm 037000}, China}

%%% Abstract. 摘要
\abstract{In this paper, we study the existence and uniqueness of pseudo $S$-asymptotically $\omega$-periodic mild solutions of class $r$ for
fractional integro-differential neutral equations. An example is presented to illustrate the application of the abstract results.}

%%% Keywords. 关键词


\maketitle

%%%%%%%%%%%%%%%%%%%%%%%%%%%%%%%%%%%%%%%%%%%%%%%%%%%%%%%%%%%%
%%%%%%%%%%%%%%%%%%%%%%%%%%%%%%%%%%%%%%%%%%%%%%%%%%%%%%%
%%% The main text. 正文部分%
%%  图表引用\cref公式引用\eqref参考文献\cite
%%%%%%%%%%%%%%%%%%%%%%%%%%%%%%%%%%%%%%%%%%%%%%%%%%%%%%%
\section{Introduction}

In the literature, several concepts were presented to study the approximate periodic function. The existence of periodic solutions to different kinds of differential equations and its various generalizations like the existence of almost periodic solutions, almost automorphic solutions, asymptotically almost periodic solutions, asymptotically almost automorphic solutions, pseudo-almost periodic solutions and pseudo-almost automorphic solutions have been extensively studied, see the monographs of Corduneanu [11], N'Gu\'{e}r\'{e}kata [29], the works [1-5, 8-10, 19, 20, 25, 27, 28, 30, 33, 34, 38, 39] and references therein.

   Asymptotic periodicity problem is an active topic of the current research in the behavior theory of solutions of differential equations. We note that many concrete systems usually are submitted to non-periodic external perturbations. In many practical situations, we can assume that these perturbations are approximately periodic in a broad sense. Recently, the $S$-asymptotically $\omega$-periodic functions were formally introduced by Henríquez et al. [24] and after that, its wild applications in functional differential equations, integro-differential equations, fractional differential equations and the existence and uniqueness of $S$-asymptotically $\omega$-periodic solutions to these equations have been well studied. See, for instance, [6,15-18,21,23,24,31]. In 2013, Pierri and Rolnik [32] introduced the concept of pseudo $S$-asymptotically $\omega$-periodic ($PSAP_{\omega}$, in short) functions, which is a natural generalization of $S$-asymptotically $\omega$-periodic functions. Also, they studied the existence and uniqueness of pseudo $S$-asymptotically $\omega$-periodic mild solutions for abstract neutral functional equations. The theory of pseudo $S$-asymptotically $\omega$-periodic functions
has been shown to have interesting applications in several branches of differential equations, and because of this fact, this theory has been attracting the attention of several mathematicians and the interest in the topic remains growing. We refer the reader to [14,22,36] and the references therein. In [6] the authors established the sufficient conditions to ensure the existence and uniqueness of pseudo $S$-asymptotically $\omega$-periodic mild solutions for hyperbolic evolution equations respectively by hyperbolic semigroup and uniformly stable semigroup and considered the $PSAP_{\omega}$-mild solutions in intermediate spaces. Xia [37] discussed the pseudo asymptotically periodic solutions of two-term time fractional differential equations
with delay. In [22] the authors established the criteria to guarantee the existence and uniqueness of pseudo $S$-asymptotically periodic solutions of second-order abstract Cauchy problems. In [14] the authors have studied qualitative properties of pseudo $S$-asymptotically $\omega-$periodic functions. Also they discussed the existence of pseudo $S$-asymptotically $\omega$-periodic mild solutions for fractional differential equations in the following form
 \begin{eqnarray}
\left \{\begin{array}{ll}\displaystyle\frac{d}{dt}v(t)=\displaystyle\int_0^t\frac{(t-s)^{\alpha-2}}{\Gamma(\alpha-1)}Av(s)ds+f(t,v(t)),~~t\geq0,\\
v(0)=v_0\in X,\end{array}
\right.
\end{eqnarray}
 where $1 < \alpha < 2$, $A : D(A)\subseteq X\rightarrow X$ is a linear densely defined operator of sectorial type on a complex Banach space $(X,\|\cdot\|)$ and $f :[0,\infty) \times X \rightarrow X$ is an appropriate function.

Fractional calculus and its applications have gained a lot of
attention due to their applications in the fields such as physics, fluid mechanics viscoelasticity, heat conduction in materials
with memory, chemistry and engineering. In
recent years, notable contributions have been made in theory and applications of fractional differential equations, one can refer to [26,40,41].

However, pseudo $S$-asymptotically $\omega$-periodic mild solutions for fractional neutral differential equations have still rarely been treated in the literature. Motivated by these facts, in this paper, we mainly investigate the existence and uniqueness of pseudo $S$-asymptotically $\omega$-periodic mild solutions of class $p$ to the system
 \begin{eqnarray}
\left \{\begin{array}{ll}\displaystyle\frac{d}{dt}D(t,u_t)=\displaystyle\int_0^t\frac{(t-s)^{\alpha-2}}{\Gamma(\alpha-1)}AD(s,u_s)ds+f(t,u_t),~~t\geq0,\\
u_0=\varphi\in \mathcal{C},\end{array}
\right.
\end{eqnarray}
 where $1 < \alpha < 2$, $D(t,\varphi)=\varphi(0)+g(t,\varphi)$, $A : D(A)\subseteq X\rightarrow X$ is a linear densely defined operator of sectorial type on a complex Banach space $X$, the history $u_t:[-r,0]\rightarrow\mathcal{C}=C([-r,0],X)$ is defined by $u_t(\theta)=u(t+\theta)$ for $\theta\in[-r,0]$ and $f,g$ are functions satisfying some additional conditions to be specified later. The convolution
integral in (1.2) is understood in the Riemann--Liouville sense.

We organize this paper as follows. In the next section, we will introduce the notions of $S$-asymptotically $\omega$-periodic and pseudo $S$-asymptotically $\omega$-periodic functions and give some basic properties. In Section 3, we obtain very general results on the existence and uniqueness of pseudo $S$-asymptotically $\omega$-periodic mild solutions for the semilinear problem (1.2) under Lipschitz type hypothesis on the nonlinearity. Finally, in Section 4, we apply our theory to concrete applications. We present an example, leading to a better understanding of the work and hence attract the attention to researchers who are entering the subject.

%%%%%%%%%%%%%
% SECTION 2 %
%%%%%%%%%%%%%

\section{Preliminaries}
\setcounter{equation}{0}\setcounter{section}{2}\indent
Assume $a$ is an arbitrary real number, $r>0$. Let $C_b([a, \infty), X),~C_b([0, \infty), \mathbb{R}^{+})$ and $\mathcal{C}=C([-r, 0], X)$ stand for the space formed by all bounded continuous functions from $[a, \infty)$ into $X$ endowed with the sup norm $\|\cdot\|_{C_b([a, \infty), X)}$, the space formed by all bounded continuous functions from $[0, \infty)$ into the positive real number set $\mathbb{R}^{+}$ endowed with the sup norm $\|\cdot\|_{C_b([0, \infty), \mathbb{R}^{+})}$ and the space of continuous functions from $[-r,0]$ into $X$ with sup norm $\|\cdot\|_{\mathcal{C}}$, respectively. The notation $\mathcal{B}(X)$ stands for the space of bounded linear operators from $X$ into
$X$ endowed with the uniform operator norm denoted $\|\cdot\|_{\mathcal{B}(X)}$.\\

A closed linear operator $A$ is said to be sectorial of type $\mu$ if there exist $0 < \theta <\frac{\pi}{2}
, M > 0$ and $\mu\in \mathbb{R}$ such that the spectrum of
$A$ is contained in the sector $\mu+S_{\theta}=\{\mu+\lambda:\lambda\in \mathbb{C}\setminus{\{0\}},|arg(-\lambda)|<\theta\}$, and $\|\lambda I-A\|^{-1}\leq\frac{M}{|\lambda-\mu|}$ for all $\lambda\notin \mu+S_{\theta}.$ (see [14])

In order to give an operator theoretical approach to system (1.2) we recall the following definition.\\[3mm]
{\bf Definition 2.1 } (See [14]). Let $A$ be a closed and linear operator with domain $D(A)$ defined on a Banach space $X$. We call $A$ the
generator of a solution operator if there exist $\mu\in \mathbb{R}$ and a strongly continuous function $S_{\alpha}:\mathbb{R}^{+}\rightarrow\mathcal{B}(X)$ such that
$\{\lambda^{\alpha}:Re\lambda>\mu\}\subseteq\rho(A)$ and $\lambda^{\alpha-1}(\lambda^{\alpha}I-A)^{-1}x=\displaystyle\int_0^{\infty}e^{-\lambda t}S_{\alpha}(t)xdt$ for all $Re\lambda>\mu,~x\in X.$
 In this case, $S_{\alpha}(t)$ is called the solution
operator generated by $A$ which satisfies $S_{\alpha}(0)=I$. We
observe that the power function $\lambda^{\alpha}$ is uniquely defined as $\lambda^{\alpha}=|\lambda^{\alpha}|e^{\alpha i\arg(\lambda)}$ with $-\pi<\arg(\lambda)<\pi$.

We note that if $A$ is a sectorial operator of type $\mu$ with $0 < \theta\leq\pi(1-\frac{\alpha}{2})$
, then $A$ is the generator of a solution operator
given by $S_{\alpha}(t):=\displaystyle\frac{1}{2\pi i}\int_{\gamma}e^{\lambda t}\lambda^{\alpha-1}(\lambda^{\alpha}I-A)^{-1}d\lambda,t>0,$ where $\gamma$ is a suitable path lying outside the sector $\mu+\Sigma_{\theta}$ (cf. [12]). In 2007,
Cuesta [12, Theorem 1] proved that if $A $ is a sectorial operator of type $\mu < 0$, for $M > 0$ and $0 < \theta\leq\pi(1-\frac{\alpha}{2})$, then there exists
$C > 0$ such that $$\|S_{\alpha}(t)\|\leq\displaystyle\frac{CM}{1+|\mu|t^{\alpha}},~~t\geq0.\eqno(2.1)$$\\
{\bf Remark 2.1.} In the rest of this paper, we always suppose that $A$ is a sectorial of type $\mu < 0$ with angle $\theta$ satisfying
$0 < \theta\leq\pi(1-\frac{\pi}{2})$, $M$ and $C$ are the constants introduced above.

We now recall some notations and properties related to $S$-asymptotically $\omega$-periodic functions.\\[3mm]
{\bf Definition 2.2 } (See [25]). A set $\Omega \subseteq \mathbb{R}$ is said to be an ergodic zero set if
$\lim\limits_{r\rightarrow\infty}\displaystyle\frac{{\rm mes} (Q_r\cap \Omega)}{{\rm mes} (Q_r)}= 0,$
where ${\rm mes}$ denotes the Lebesgue measure and $Q_r=[-r,r].$\\[3mm]
{\bf Definition 2.3 } (See [19]). The Bochner transform $f^b
(t, s), t \in \mathbb{R},s \in[0, 1]$, of a functions $f :\mathbb{R}\rightarrow X$ is defined by
$f^b
(t,s):=f(t+s).$\\[3mm]
{\bf Definition 2.4 } (See [19]). Let $p\in[1,+\infty)$, the space $BS^p(X)$ of all Stepanov bounded functions, with the exponent $p$, consists of all
measurable functions $f :\mathbb{R}\rightarrow X$ such that $f^b
\in L^{\infty}(\mathbb{R},L^p(0,1;X))$. This is a Banach space with the norm
$$\|f\|_{S^{p}}=\|f^{b}\|_{L^{\infty}(\mathbb{R},L^p)}=\sup\limits_{t\in \mathbb{R}}\bigg(\int_{t}^{t+1}\|f(\tau)\|^pd\tau\bigg)^{\frac{1}{p}}.$$
{\bf Definition 2.5 } (See [24]).  A function $f\in C_b([0,\infty),X)$  is said to be $S$-asymptotically $\omega$-periodic if $\lim\limits_{t\rightarrow\infty}\|f(t+\omega)-f(t)\|=0$. In this case, we say that $\omega$ is an asymptotic period of $f$.\\

We use the notation $SAP_\omega(X)$ (respectively, $AP_\omega(X))$ to represent the subspace of $C_b([0,+\infty),X)$ formed by all $S$-asymptotically $\omega$-periodic functions (respectively, asymptotically $\omega$-periodic). We note that $SAP_\omega(X)$ and $AP_\omega(X)$ endowed with the norm of uniform convergence are Banach spaces.\\[3mm]
{\bf Definition 2.6 } (See [14]).  A function $f\in C_b([0,\infty), X)$ (respectively, $f\in C_b(\mathbb{R}, X)$) is said to be pseudo $S$-asymptotically $\omega$-periodic if $\displaystyle\lim\limits_{t\rightarrow\infty}\frac{1}{t}\int_{0}^{t}\|f(s+\omega)-f(s)\|ds=0$ (respectively, $\displaystyle\lim\limits_{t\rightarrow\infty}\frac{1}{2t}\int_{-t}^{t}\|f(s+\omega)-f(s)\|ds=0).$
In this case, we say that $\omega$ is an asymptotic period of $f$.\\

We use the notation $PSAP_\omega(X)$ to represent the subspace of $C_b([0,\infty),X)$ formed by all pseudo $S$-asymptotically $\omega$-periodic functions. We observe that $PSAP_\omega(X)$ endowed with the norm of uniform convergence is a Banach space and $AP_{\omega}(X)\hookrightarrow SAP_{\omega}(X)\hookrightarrow PSAP_{\omega}(X)$.\\[3mm]
{\bf Remark 2.2 }(See [14]). We observe that $u \in PSAP_{\omega}(X)$ if and only if for each $\varepsilon>0$, the set
$C_{\varepsilon}=\{t\in [0,\infty):\|u(t+\omega)-u(t)\|\geq\varepsilon\}$
is an ergodic zero set.\\[3mm]
{\bf Definition 2.7 } (See [32]). Let $r > 0$ and $u \in PSAP_{\omega}(X)$. We say that $u$ is pseudo $S$-asymptotically $\omega$-periodic of class $r$ if
 $$\lim\limits_{T\rightarrow\infty}\frac{1}{T}\int_{r}^{T}\sup\limits_{\tau\in [s-r,s]}\|u(\tau+\omega)-u(\tau)\|ds=0.\eqno(2.2)$$

 We denote by $PSAP_{\omega,r}(X)$ the subspace of $C_b([0, \infty), X)$ formed by all the pseudo $S$-asymptotically $\omega$-periodic functions of class $r$. We note that $PSAP_{\omega,r}(X)$ endowed with the norm of uniform convergence is a Banach space.\\[3mm]
 {\bf Remark 2.3 }(See [37]). Assume $u \in C_b([0,+\infty),X)$, then $u \in PSAP_{\omega,r}(X)$ if and only for $\forall~\varepsilon>0,$ \\
 $~~~~~~~~~\lim\limits_{T\rightarrow\infty}\displaystyle\frac{ 1}{T}{\rm mes}(M_{T,\varepsilon}(u))=0,$\\[3mm]
  where ${\rm mes}$ denotes the Lebesgue measure and\\[2mm]
$~~~~~~~~~M_{T,\varepsilon}(u)=\bigg\{t\in [r,T]:\sup\limits_{\tau\in[t-r,t]}\|u(\tau+\omega)-u(\tau)\|\geq\varepsilon\bigg\}.$\\[3mm]
 {\bf Lemma 2.1 } (See [37]). \emph{Assume $r\geq0,$ $r_1>0,~r_2>0,$ then\\
(i) $PSAP_{\omega,r}(X)\subseteq PSAP_{\omega}(X).$\\
(ii) $PSAP_{\omega,r}(X)$ is a closed subspace of $C_b([0,+\infty), X).$\\
(iii) $PSAP_{\omega,r_1}(X)=PSAP_{\omega,r_2}(X)$}.\\[3mm]
{\bf Lemma 2.2.} \emph{The space $PSAP_{\omega,r}(X)$ is translation invariant on $\mathbb{R}^{+}$.}\\[3mm]
\emph{ Proof.} Let $u\in PSAP_{\omega,r}(X)$. For $s\in\mathbb{ R}^{+},$ we have the estimate
\begin{eqnarray*}
&&\frac{1}{t}\int_r^t\sup\limits_{\rho\in[\xi-r,\xi]}\|u(\rho+\omega+s)-u(\rho+s)\|d\xi\\
&&=(1+\frac{s}{t})\bigg(\frac{1}{t+s}\int_{r+s}^{t+s}\sup\limits_{\rho\in[\xi-r,\xi]}\|u(\rho+\omega)-u(\rho)\|d\xi\bigg)\\
&&\leq(1+\frac{s}{t})\bigg(\frac{1}{t+s}\int_r^{t+s}\sup\limits_{\rho\in[\xi-r,\xi]}\|u(\rho+\omega)-u(\rho)\|d\xi\bigg).
\end{eqnarray*}
Thus we show that the function $\rho\rightarrow u(\rho+s)$ belongs to $PSAP_{\omega,r}(X)$.$\quad\Box$\\[3mm]
{\bf Definition 2.8 }(See [32]). We say that a function $F \in C([0, \infty)\times \mathcal{C}, X)$ is uniformly $(\mathcal{C},X)$
pseudo $S$-asymptotically $\omega$-periodic of class $r$ if
$$\lim\limits_{T\rightarrow\infty}\frac{1}{T}\int_{r}^{T}\sup\limits_{\tau\in [s-r,s]}\sup\limits_{\|x\|_{\mathcal{C}}\leq L}\|F(\tau+\omega,x)-F(\tau,x)\|ds=0,\eqno(2.3)$$
for all $L>0$. Denote by $PSAP_{\omega,r}(\mathcal{C},X)$ the set formed by functions of this type.\\[3mm]
{\bf Lemma 2.3 }(See [32]). \emph{Let $u\in C_b([-r, \infty), X)$ and assume that $u\mid_{[0,\infty)} \in PSAP_{\omega,r}(X)$. Then the
function $s \rightarrow u_s$ belongs to $PSAP_{\omega,r}(\mathcal{C})$.}\\[3mm]
{\bf Lemma 2.4.} \emph{Assume that $F\in PSAP_{\omega,r}(\mathcal{C},X)$ and there exists $L_F\in C_b([0, \infty), \mathbb{R}^{+})$ such that
$\|F(t,\psi_{1}) - F(t,\psi_{2})\|\leq L_F(t)\| \psi_{1} - \psi_{2} \|_{\mathcal{C}}$ for all $(t,\psi_{i}) \in [0,\infty) \times \mathcal{C}$. If $u \in
C_b([-r, \infty), X)$ and $u\mid_{[0,\infty)}\in PSAP_{\omega,r}(X)$, then the function $s \rightarrow F(s, u_s)$
belongs to $PSAP_{\omega,r}(X)$.}\\[3mm]
\emph{ Proof.} Let $Q=\|u\|_{C_b([-r,\infty),X)}$. Since the function
$s\mapsto u_s$ belongs to $PSAP_{\omega,r}(\mathcal{C})$ (see Lemma 2.3), then for $\forall~\varepsilon > 0$ there exists $T_{\varepsilon} > 0$ such that for each $T_{\varepsilon}>T,$ we have
$$\|L_F\|_{C_b([0,\infty),\mathbb{R}^{+})}\frac{1}{T}\int_{r}^{T}\sup\limits_{t\in[s-r,s]}\|u_{t+\omega}-u_{t}\|_{\mathcal{C}}ds\leq\varepsilon.$$
Moreover, by $F\in PSAP_{\omega,r}(\mathcal{C},X)$, we have
$$\frac{1}{T}\int_{r}^{T}\sup\limits_{t\in[s-r,s]}\sup\limits_{\|x\|_{\mathcal{C}}\leq Q}\|F(t+\omega,x)-F(t,x)\|ds\leq\varepsilon.$$
Thus for all $T\geq T_{\varepsilon},$ we can obtain
\begin{eqnarray*}
&&\frac{1}{T}\int_{r}^{T}\sup\limits_{t\in[s-r,s]}\|F(t+\omega,u_{t+\omega})-F(t,u_{t})\|ds\\
&&~\leq\frac{1}{T}\int_{r}^{T}\sup\limits_{t\in[s-r,s]}\|F(t+\omega,u_{t+\omega})-F(t,u_{t+\omega})\|ds\\
&&~~~+\frac{1}{T}\int_{r}^{T}\sup\limits_{t\in[s-r,s]}\|F(t,u_{t+\omega})-F(t,u_{t})\|ds\\
&&~\leq\frac{1}{T}\int_{r}^{T}\sup\limits_{t\in[s-r,s]}\sup\limits_{\|x\|_{\mathcal{C}}\leq Q}\|F(t+\omega,x)-F(t,x)\|ds\\
&&~~~~+\|L_F\|_{C_b([0,\infty);\mathbb{R}^{+})}\frac{1}{T}\int_{r}^{T}\sup\limits_{t\in[s-r,s]}\|u_{t+\omega}-u_{t}\|_{\mathcal{C}}ds\\
&&~\leq2\varepsilon,
\end{eqnarray*}
which claims the assertion.$\quad\Box$\\[3mm]
{\bf Lemma 2.5.} \emph{Let $u \in C_b([-r, \infty), X)$ and $\kappa : [-r,\infty)\rightarrow X$ be the function defined by $\kappa(t) =0$ for $t\in[-r,0]$ and $\kappa(t) =
\displaystyle\int_0^t S_{\alpha}(t-s)u(s)ds$ for $t \geq 0$. If $u\mid_{[0,\infty)} \in PSAP_{\omega,r}(X)$, then $\kappa \in PSAP_{\omega,r}(X)$.}\\[3mm]
\emph{ Proof.} We have the estimate
$$\int_0^{t}\|S_{\alpha}(t-s)u(s)\|ds\leq\frac{CM|\mu|^{\frac{-1}{\alpha}}\pi}{\alpha\sin(\frac{\pi}{\alpha})}\|u\|_{C_b([0,\infty),X)},\eqno(2.4)$$
which shows that $\kappa\in C_b([-r,\infty),X)$. Now we shall show that $$\lim\limits_{ T\rightarrow\infty}\displaystyle\frac{1}{T}\int_r^T\sup\limits_{t\in[\xi-r,\xi]}\|\kappa(t+\omega)-\kappa(t)\|d\xi= 0.$$

We first choose positive constant $q$ such that $T>q>r>0$, then we can obtain
\begin{eqnarray*}
&&\frac{1}{T}\int_r^T\sup\limits_{t\in[\xi-r,\xi]}\|\kappa(t+\omega)-\kappa(t)\|d\xi\\
&&~\leq\frac{1}{T}\int_r^T\sup\limits_{t\in[\xi-r,\xi]}\int_0^t\|S_{\alpha}(t-s)[u(s+\omega)-u(s)]\|dsd\xi\\
&&~~~~+\frac{1}{T}\int_r^T\sup\limits_{t\in[\xi-r,\xi]}\int_{0}^{\omega}\|S_{\alpha}(t+\omega-s)u(s)\|dsd\xi\\
&&~\leq\frac{1}{T}\int_r^q\sup\limits_{t\in[\xi-r,\xi]}\int_0^t\|S_{\alpha}(t-s)[u(s+\omega)-u(s)]\|dsd\xi\\
&&~~~~+\frac{1}{T}\int_q^T\sup\limits_{t\in[\xi-r,\xi]}\int_0^t\|S_{\alpha}(t-s)[u(s+\omega)-u(s)]\|dsd\xi\\
&&~~~~+\frac{1}{T}\int_r^T\sup\limits_{t\in[\xi-r,\xi]}\int_0^{\omega}\|S_{\alpha}(t+\omega-s)u(s)\|dsd\xi:=\sum\limits_{i=1}^{3}K_{i}(T).
\end{eqnarray*}

Next, we will estimate the terms $K_i(T),1\leq i\leq3,$ separately.

For the term $K_1(T)$, we aim to prove that $\lim\limits_{T\rightarrow\infty}K_1(T)=0.$ Since $u\in C_{b}([-r,\infty),X)$, then we can get
\begin{eqnarray*}
K_{1}(T)&=&\frac{1}{T}\int_r^q\sup\limits_{t\in[\xi-r,\xi]}\int_0^{t}\|S_{\alpha}(t-s)[u(s+\omega)-u(s)]\|dsd\xi\\
&\leq&\frac{1}{T}\int_r^q\sup\limits_{t\in[\xi-r,\xi]}\int_0^{t}\frac{CM}{1+|\mu|(t-s)^{\alpha}}[u(s+\omega)-u(s)]\|dsd\xi\\
&\leq&2CM\|u\|_{C_b([-r,\infty),X)}\frac{q(q-r)}{T}\rightarrow0,  \mbox{\quad as }T\rightarrow\infty.
\end{eqnarray*}

To verify $\lim\limits_{T\rightarrow\infty}K_2(T)=0,$ we first decompose $K_{2}(T)$ as follows
\begin{eqnarray*}
K_{2}(T)&=&\frac{1}{T}\int_q^{T}\sup\limits_{t\in[\xi-r,\xi]}\int_0^t\|S_{\alpha}(t-s)[u(s+\omega)-u(s)]\|dsd\xi\\
&\leq&\frac{1}{T}\int_q^{T}\sup\limits_{t\in[\xi-r,\xi]}\int_0^{t-2}\|S_{\alpha}(t-s)[u(s+\omega)-u(s)]\|dsd\xi\\
&&+\frac{1}{T}\int_q^{T}\sup\limits_{t\in[\xi-r,\xi]}\int_{t-2}^t\|S_{\alpha}(t-s)[u(s+\omega)-u(s)]\|dsd\xi\\
&=&\sum\limits_{i=1}^{2}K_{2}^{i}(T),
\end{eqnarray*}
where $$ K_{2}^{1}(T)=\displaystyle\frac{1}{T}\int_q^{T}\sup\limits_{t\in[\xi-r,\xi]}\displaystyle\int_0^{t-2}\|S_{\alpha}(t-s)[u(s+\omega)-u(s)]\|dsd\xi,$$ and $$K_{2}^{2}(T)=\displaystyle\frac{1}{T}\int_q^{T}\sup\limits_{t\in[\xi-r,\xi]}\displaystyle\int_{t-2}^t\|S_{\alpha}(t-s)[u(s+\omega)-u(s)]\|dsd\xi.$$

For the term $K_{2}^{1}(T)$, we have
\begin{eqnarray*}
K_{2}^{1}(T)&=&\frac{1}{T}\int_q^{T}\sup\limits_{t\in[\xi-r,\xi]}\int_0^{t-2}\|S_{\alpha}(t-s)[u(s+\omega)-u(s)]\|dsd\xi\\
&\leq&\frac{1}{T}\int_q^{T}\sup\limits_{t\in[\xi-r,\xi]}\int_{0}^{t-2}\frac{CM}{1+|\mu| (t-s)^{\alpha}}\|u(s+\omega)-u(s)\|dsd\xi\\
&\leq&\frac{1}{T}\int_{0}^{T}\int_{0}^{\xi}\frac{CM}{1+|\mu|(\xi-r-s)^{\alpha}}\|u(s+\omega)-u(s)\| dsd\xi\\
&\leq&\frac{1}{T}\int_{0}^{T}\int_{s}^{T}\frac{CM}{1+|\mu|(\xi-r-s)^{\alpha}}\|u(s+\omega)-u(s)\|d\xi ds\\
&\leq&\frac{1}{T}\frac{CM}{|\mu|}\int_{0}^{T}[\frac{(T-r-s)^{1-\alpha}}{|1-\alpha|}+\frac{r^{1-\alpha}}{|1-\alpha|}]\|u(s+\omega)-u(s)\| ds\\
&\leq&\frac{1}{T}\frac{CM}{|\mu|}(\int_{0}^{r}+\int_{r}^{T})[\frac{(T-r-s)^{1-\alpha}}{|1-\alpha|}+\frac{r^{1-\alpha}}{|1-\alpha|}]\|u(s+\omega)-u(s)\| ds\\
&\leq&\frac{2}{T}\frac{CM}{|\mu|}\bigg[\frac{2(T-2r)^{2-\alpha}}{(\alpha-1)(2-\alpha)}+\frac{(T-r)^{2-\alpha}}{(\alpha-1)(2-\alpha)}+\frac{r^{2-\alpha}}{|1-\alpha|}+\frac{r^{2-\alpha}}{(\alpha-1)(2-\alpha)}\bigg]\\
&&\times\|u\|_{C_b([0,\infty),X)}+\frac{CM}{|\mu|}\frac{r^{1-\alpha}}{|1-\alpha|}\frac{1}{T}\int_{r}^{T}\|u(s+\omega)-u(s)\|ds.
\end{eqnarray*}
From the estimate $\displaystyle\frac{1}{T}\int_{r}^{T}\|u(s+\omega)-u(s)\|ds\leq\frac{1}{T}\int_{r}^{T}\sup\limits_{s\in[\tau-r,\tau]}\|u(s+r+\omega)-u(s+r)\|d\tau$ and $PSAP_{\omega,r}(X)$ is translation invariant on $\mathbb{R}^{+}$, we can infer that $$\lim\limits_{T\rightarrow\infty}\displaystyle\frac{1}{T}\int_{r}^{T}\|u(s+\omega)-u(s)\|ds=0.\eqno(2.5)$$ Also, since $1<\alpha<2$, then we can easily get $\lim\limits_{T\rightarrow\infty}K_{2}^{1}(T)=0.$

Now the term $K_{2}^{2}(T)$ can be estimated as,
$$
\begin{array}{rl}
K_{2}^{2}(T)&=\displaystyle\frac{1}{T}\int_q^{T}\sup\limits_{t\in[\xi-r,\xi]}\int_{t-2}^t\|S_{\alpha}(t-s)[u(s+\omega)-u(s)]\|dsd\xi\\[4mm]
&\leq\displaystyle\frac{CM}{|\mu|}\displaystyle\frac{1}{T}\int_{0}^{T} \int_{\xi-r-2}^{\xi}\frac{1}{(\xi-r-s)^{\alpha}}\|u(s+\omega)-u(s)\|ds d\xi\\[4mm]
&=\displaystyle\frac{CM}{|\mu|}\frac{1}{T}\bigg[\int_{-r-2}^{0} \int_{0}^{s+r+2}\frac{1}{(\xi-r-s)^{\alpha}}\| u(s+\omega)-u(s)\|d\xi ds\\[4mm]
&~+\displaystyle\int_{0}^{T-r-2} \int_{s}^{s+r+2}\frac{1}{(\xi-r-s)^{\alpha}}\| u(s+\omega)-u(s)\|d\xi ds\\[4mm]
&~+\displaystyle\int_{T-r-2}^{T}\int_{s}^{T} \frac{1}{(\xi-r-s)^{\alpha}}\| u(s+\omega)-u(s)\|d\xi ds\bigg]\\[4mm]
&=\displaystyle\frac{CM}{|\mu|}\bigg[\frac{1}{T}\int_{-r-2}^{0} [\frac{2^{1-\alpha}}{1-\alpha}-\frac{(-r-s)^{1-\alpha}}{1-\alpha}] \|u(s+\omega)-u(s)\|ds\\[4mm]
&~+\displaystyle\frac{1}{T}\int_{0}^{T-r-2}[\frac{2^{1-\alpha}}{1-\alpha}-\frac{(-r)^{1-\alpha}}{1-\alpha}]\|u(s+\omega)-u(s)\|ds\\[4mm]
&~+\displaystyle\frac{1}{T}\int_{T-r-2}^{T}[\frac{(T-r-s)^{1-\alpha}}{1-\alpha}-\frac{(-r)^{1-\alpha}}{1-\alpha}]\|u(s+\omega)-u(s)\|ds\bigg]\\[4mm]
&=\displaystyle\frac{CM}{|\mu|}\frac{1}{T}\int_{-r-2}^{0} \frac{2^{1-\alpha}}{1-\alpha}\|u(s+\omega)-u(s)\|ds\\
&~+\displaystyle\frac{CM}{|\mu|}\frac{1}{T}\displaystyle\int_{0}^{T-r-2}[\frac{2^{1-\alpha}}{1-\alpha}-\frac{(-r)^{1-\alpha}}{1-\alpha}]\|u(s+\omega)-u(s)\|ds\\
&~+\displaystyle\frac{CM}{|\mu|}\frac{1}{T}\displaystyle\int_{T-r-2}^{T}\frac{(-r)^{1-\alpha}}{\alpha-1}\|u(s+\omega)-u(s)\|ds\\[4mm]
&~+\displaystyle\frac{CM}{|\mu|}\frac{1}{T}\int_{-r-2}^{0} \frac{(-r-s)^{1-\alpha}}{\alpha-1} \|u(s+\omega)-u(s)\|ds\\
&~+\displaystyle\frac{CM}{|\mu|}\displaystyle\frac{1}{T}\int_{T-r-2}^{T}\frac{(T-r-s)^{1-\alpha}}{1-\alpha}\|u(s+\omega)-u(s)\|ds\\
&=\sum\limits_{i=1}^{5}I_{i}.
\end{array}
$$
Now we estimate $I_{i}$ separately.
\begin{eqnarray*}
I_{1}&=&\displaystyle\frac{CM}{|\mu|}\frac{1}{T}\int_{-r-2}^{0} \frac{2^{1-\alpha}}{1-\alpha}\|u(s+\omega)-u(s)\|ds\\
&\leq&\displaystyle\frac{CM}{|\mu|}\frac{2(r+2)}{T}\frac{2^{1-\alpha}}{1-\alpha}\|u\|_{C_b([0,\infty),X)}\rightarrow0, \mbox{\quad as }T\rightarrow0.
\end{eqnarray*}
For $I_{2},$ by a standard calculation and $u\mid_{[0,\infty)} \in PSAP_{\omega,r}(X)$, we have
\begin{eqnarray*}
I_{2}&=&\displaystyle\frac{CM}{|\mu|}\frac{1}{T}\displaystyle\int_{0}^{T-r-2}[\frac{2^{1-\alpha}}{1-\alpha}-\frac{(-r)^{1-\alpha}}{1-\alpha}]\|u(s+\omega)-u(s)\|ds\\
&\leq&\displaystyle\frac{CM}{|\mu|}\frac{1}{T}\displaystyle\int_{0}^{T-r}[\frac{2^{1-\alpha}}{1-\alpha}-\frac{(-r)^{1-\alpha}}{1-\alpha}]\|u(s+\omega)-u(s)\|ds\\
&\leq&\displaystyle\frac{CM}{|\mu|}[\frac{2^{1-\alpha}}{1-\alpha}-\frac{(-r)^{1-\alpha}}{1-\alpha}]\frac{2}{T}\displaystyle[\int_{0}^{r}\|u(s+\omega)-u(s)\|+\int_{r}^{T}\|u(s+\omega)-u(s)\|ds]\\
&\leq&\displaystyle\frac{CM}{|\mu|}[\frac{2^{1-\alpha}}{1-\alpha}-\frac{(-r)^{1-\alpha}}{1-\alpha}]\frac{2}{T}\displaystyle[2r\|u\|_{C_b([0,\infty),X)}+\int_{r}^{T}\|u(s+\omega)-u(s)\|ds]\rightarrow0,\mbox{\quad as }T\rightarrow0.
\end{eqnarray*}
For $I_{3},$ since $\|u(s+\omega)-u(s)\|\leq2\|u\|_{C_b([0,\infty),X)}$, thus we can deduce
\begin{eqnarray*}
I_{3}&=&\displaystyle\frac{CM}{|\mu|}\frac{1}{T}\displaystyle\int_{T-r-2}^{T}\frac{(-r)^{\alpha-1}}{1-\alpha}\|u(s+\omega)-u(s)\|ds\\
&\leq&\frac{CM}{|\mu|}\frac{2}{T}(r+2)\frac{(-r)^{1-\alpha}}{\alpha-1}\|u\|_{C_b([0,\infty),X)}\rightarrow0,\mbox{\quad as }T\rightarrow0.
\end{eqnarray*}
For $I_{4}$ and $I_{5}$, we can derive
\begin{eqnarray*}
I_{4}&=&\displaystyle\frac{CM}{|\mu|}\frac{1}{T}\int_{-r-2}^{0} \frac{(-r-s)^{1-\alpha}}{\alpha-1} \|u(s+\omega)-u(s)\|ds\\
&=&\displaystyle\frac{CM}{|\mu|}\frac{2}{T}[\frac{(-r)^{(2-\alpha)(1-\alpha)}}{(2-\alpha)(\alpha-1)}-\frac{2^{(2-\alpha)(1-\alpha)}}{(2-\alpha)(\alpha-1)}]\|u\|_{C_b([0,\infty),X)},\\
I_{5}&=&\displaystyle\frac{CM}{|\mu|}\displaystyle\frac{1}{T}\int_{T-r-2}^{T}\frac{(T-r-s)^{1-\alpha}}{1-\alpha}\|u(s+\omega)-u(s)\|ds\\
&=&\displaystyle\frac{CM}{|\mu|}\frac{2}{T}[\frac{(-r)^{(2-\alpha)(1-\alpha)}}{(2-\alpha)(1-\alpha)}-\frac{2^{(2-\alpha)(1-\alpha)}}{(2-\alpha)(1-\alpha)}]\|u\|_{C_b([0,\infty),X)}.
\end{eqnarray*}
We can get $I_{4}+I_{5}=0.$

Combing with the estimates above, expression (2.5) and $1<\alpha<2$, we can infer that $\lim\limits_{T\rightarrow\infty}K_{2}^2(T)=0.$

Finally,
\begin{eqnarray*}
K_{3}(T)&=&\frac{1}{T}\int_r^T\sup\limits_{t\in[\xi-r,\xi]}\int_{0}^{\omega}\|S_{\alpha}(t+\omega-s)u(s)\|dsd\xi\\
&\leq&\frac{1}{T}\int_r^T\frac{CM}{1+|\mu|(\xi-r)^{\alpha}}\int_{0}^{\omega}\|u(s)\|dsd\xi\\
&\leq&\frac{1}{T}\omega\|u\|_{C_b([0,\infty),X)}\int_r^T\frac{CM}{1+|\mu|(\xi-r)^{\alpha}}d\xi\\
&\leq&\frac{1}{T}\omega\|u\|_{C_b([0,\infty),X)}\frac{CM|\mu|^{\frac{-1}{\alpha}}\pi}{\alpha\sin(\frac{\pi}{\alpha})}\rightarrow0,  \mbox{\quad as }T\rightarrow\infty.
\end{eqnarray*}

Consequently, we prove that $\kappa \in PSAP_{\omega,r}(X)$.$\quad\Box$

%%%%%%%%%%%%%
% SECTION 3%
%%%%%%%%%%%%%

\section{Existence and uniqueness of pseudo $S$-asymptotically $\omega$-periodic solutions}
\setcounter{equation}{0}\setcounter{section}{3}\indent
This section is mainly concerned with the existence and uniqueness results of pseudo $S$-asymptotically $\omega$-periodic mild solutions of class $r$. \\[3mm]
{\bf Definition 3.1.} A function $u \in C_b([-r, \infty), X)$ is said to be a mild solution of system (1.2)
if $u_0 = \varphi$ and

$$
u(t)=S_{\alpha}(t)[\varphi(0)-g(0,\varphi)]+g(t,u_t)+\int_{0}^t S_{\alpha}(t-s)f(s,u_s)ds,~~~\forall~t\geq0.
$$
To establish our results, we require the following assumptions:

(H$_1$) The functions $g,f \in PSAP_{\omega,r}(\mathcal{C}, X)$ and $f(\cdot,0),g(\cdot,0)$ belong to $C_{b}([0,\infty),X)$;

(H$_2$) There exist positive constants $L_g,L_f$ such that for any
$\psi, \varphi\in \mathcal{C}$ and all $t \in [0,\infty)$, $$\|g(t,\psi)-g(t,\varphi)\|\leq L_g\|\psi-\varphi\|_{\mathcal{C}},~\|f(t,\psi)-f(t,\varphi)\|\leq L_f\|\psi-\varphi\|_{\mathcal{C}};$$

(H$_3$) There exist $L_g(t),L_f(t)\in C_b([0,\infty),\mathbb{R}^{+})$ such that for any
$\psi, \varphi\in \mathcal{C}$ and all $t \in [0,\infty)$, $$\|g(t,\psi)-g(t,\varphi)\|\leq L_g(t)\|\psi-\varphi\|_{\mathcal{C}},\|f(t,\psi)-f(t,\varphi)\|\leq L_f(t)\|\psi-\varphi\|_{\mathcal{C}};$$

(H$_4$) There exist $L_g(t)\in  C_b([0,\infty),\mathbb{R}^{+})$ and $L_f(t)\in C_{b}([0,\infty),\mathbb{R}^{+})\bigcap L_{loc}^1([0,+\infty),\mathbb{R}^{+})$ such that for any
$\psi, \varphi\in \mathcal{C}$ and all $t \in [0,\infty)$, $$\|g(t,\psi)-g(t,\varphi)\|\leq L_g(t)\|\psi-\varphi\|_{\mathcal{C}},\|f(t,\psi)-f(t,\varphi)\|\leq L_f(t)\|\psi-\varphi\|_{\mathcal{C}};$$

(H$_5$) There exist $L_g(t)\in  C_b([0,\infty),\mathbb{R}^{+})$ and $L_f(t)\in BS^p(\mathbb{R}^{+},\mathbb{R}^{+})\cap L_{loc}^1([0,+\infty),\mathbb{R}^{+})$ such that for any
$\psi, \varphi\in \mathcal{C}$ and all $t \in [0,\infty)$, $$\|g(t,\psi)-g(t,\varphi)\|\leq L_g(t)\|\psi-\varphi\|_{\mathcal{C}},\|f(t,\psi)-f(t,\varphi)\|\leq L_f(t)\|\psi-\varphi\|_{\mathcal{C}}.$$

We are now in a position to establish our first existence theorem.\\[3mm]
{\bf Theorem 3.1.} \emph{Assume that hypotheses} (H$_1$) \emph{and} (H$_3$) \emph{hold, and  $$\|L_g\|_{C_b([0,\infty),\mathbb{R}^{+})}+\|L_f\|_{C_b([0,\infty),\mathbb{R}^{+})}\frac{CM|\mu|^{\frac{-1}{\alpha}}\pi}{\alpha\sin(\frac{\pi}{\alpha})}<1,\eqno(3.1)$$ then system $(1.2)$ has a unique mild solution $u\in PSAP_{\omega,r}(X)$}.\\[3mm]
\emph{ Proof.} Let $\mathcal{B}=\{u:[-r,\infty)\rightarrow X~|~u_0=\varphi,u|_{[0,\infty)}\in PSAP_{\omega,r}(X)\}$ endowed with the
metric $d(u,z)=\|u-z\|_{C_b([-r,\infty),X)}$ and let $\Phi_{\alpha}$ be the map defined by
$$
\begin{array}{rl}
\Phi_{\alpha} u(t)=S_{\alpha}(t)[\varphi(0)-g(0,\varphi)]+g(t,u_t)+\displaystyle\int_{0}^t S_{\alpha}(t-s)f(s,u_s)ds.~~~\forall~t\geq0.
\end{array}\eqno(3.2)$$
 We next prove that $\Phi_{\alpha}$ is a contraction on $\mathcal{B}$.

Assume $u\in\mathcal{B},$ we now show that $\Phi_{\alpha}u\in\mathcal{B}$.

From (H$_3$) we know that $f,g$ satisfy the Lipschitz condition. Moreover, by a straightforward computation, we obtain the following estimate
$$
\begin{array}{rl}
\|\Phi_{\alpha}  u(t)\|&\leq \|S_{\alpha}(t)\|[\|\varphi(0)\|+\|g(0,\varphi)\|]+\bigg[\|g(t, u_{t})-g(t,0)\|+\|g(\cdot,0)\|_{C_b([0,\infty),X)}\bigg]\\[3mm]
&~~+\displaystyle\int_0^tS_{\alpha}(t-s)\bigg[\|f(s, u_{s})-f(s,0)\|+\|f(\cdot,0)\|_{C_b([0,\infty),X)}\bigg]ds\\
&\leq CM[\|\varphi(0)\|+\|g(0,\varphi)\|]+\bigg[\|L_g\|_{C_b([0,\infty),\mathbb{R}^{+})}\|u\|_{C_b([-r,\infty),X)}+\|g(\cdot,0)\|_{C_b([0,\infty),X)}\bigg]\\[3mm]
&~~+\displaystyle\frac{CM|\mu|^{\frac{-1}{\alpha}}\pi}{\alpha\sin(\frac{\pi}{\alpha})}\bigg[\|L_f\|_{C_b([0,\infty),\mathbb{R}^{+})} \|u\|_{C_b([-r,\infty),X)}+\|f(\cdot,0)\|_{C_b([0,\infty),X)}\bigg].
\end{array}\eqno(3.3)$$

On the other hand, it is obvious that
$$
\begin{array}{rl}
&\displaystyle\frac{1}{T}\int_0^T\sup\limits_{\tau\in[s-r,s]}\bigg\|S_{\alpha}(\tau+\omega)[\varphi(0)-g(0,\varphi)]-S_{\alpha}(\tau)[\varphi(0)-g(0,\varphi)]\bigg\|ds\\
&~\leq\displaystyle\frac{2CM}{T}\int_{0}^{T}\frac{1}{1+|\mu|\tau^{\alpha}}d\tau[\|\varphi(0)\|+\|g(0,\varphi)\|]\\
&~\leq\displaystyle\frac{2CM|\mu|^{\frac{-1}{\alpha}}\pi}{\alpha\sin(\frac{\pi}{\alpha})}\times\displaystyle\frac{1}{T}[\|\varphi(0)\|+\|g(0,\varphi)\|]\rightarrow0,  \mbox{\quad as }T\rightarrow\infty.
\end{array}\eqno(3.4)$$

Hence, we infer that $S_{\alpha}(\cdot)[\varphi(0)-g(0,\varphi)]\in PSAP_{\omega,r}(X).$
Using now the estimates (3.3), (3.4) and Lemmas 2.4-2.5, we can claim that $\Phi_{\alpha}u\in \mathcal{B}$ whenever $u\in \mathcal{B}$.

Furthermore, for all $u,v\in \mathcal{B}$ and $t\geq0$, we can deduce that
$$
\begin{array}{rl}
&\|\Phi_{\alpha}u(t)-\Phi_{\alpha}v(t)\|\\
&\leq L_g(t)\|u_t-v_t\|_{\mathcal{C}}+ CM\displaystyle\int_{0}^{t}\displaystyle\frac{L_f(s)}{1+|\mu|(t-s)^{\alpha}}ds \|u_s-v_s\|_{\mathcal{C}}\\
&\leq L_g(t)\|u_t-v_t\|_{\mathcal{C}}+\|L_f\|_{C_b([0,\infty),\mathbb{R}^{+})}\displaystyle\frac{CM|\mu|^{\frac{-1}{\alpha}}\pi}{\alpha\sin(\frac{\pi}{\alpha})}\|u_s-v_s\|_{\mathcal{C}}\\
&\leq\bigg[\|L_g\|_{C_b([0,\infty),\mathbb{R}^{+})}+\|L_f\|_{C_b([0,\infty),\mathbb{R}^{+})}\displaystyle\frac{CM|\mu|^{\frac{-1}{\alpha}}\pi}{\alpha\sin(\frac{\pi}{\alpha})}\bigg]\|u-v\|_{C_b([-r,\infty),X)},
\end{array}$$
which follows from (3.1) that $\Phi_{\alpha}$ is a contraction on $\mathcal{B}$. Therefore we can affirm that $\Phi_{\alpha}$ has a unique fixed point $u\in PSAP_{\omega,r}(X)$ in $\mathcal{B}$, which is the mild solution of (1.2). The proof is complete.$\quad\Box$

The next result is an immediate consequence of Theorem 3.1.\\[3mm]
{\bf Corollary 3.2.} \emph{Assume that hypotheses} (H$_1$) \emph{and} (H$_2$) \emph{hold, and
$$L_g+\displaystyle\frac{CML_f|\mu|^{\frac{-1}{\alpha}}\pi}{\alpha\sin(\frac{\pi}{\alpha})}<1,$$
then system $(1.2)$ has a unique mild solution $u\in PSAP_{\omega,r}(X)$}.

A similar result can be established when $f$ satisfies a local Lipschitz condition.\\[3mm]
{\bf Theorem 3.3.} \emph{Assume that hypotheses} (H$_1$) \emph{and} (H$_4$) \emph{hold, and $$\sup\limits_{t\geq0}W_f(t)+\|L_{g}\|_{C_b([0,\infty),\mathbb{R}^{+})}<1,\eqno(3.5)$$
where $W_f(t)=\displaystyle\int_0^t\frac{L_f(s)}{1+|\mu|(t-s)^{\alpha}}ds$. Then system (1.2) has a unique mild solution $u\in PSAP_{\omega,r}(X)$.}\\[3mm]
\emph{ Proof.} We still define the map $\Phi_{\alpha}$ and the set $\mathcal{B}$ as in Theorem 3.1. Next, we prove that $\Phi_{\alpha}$ is a contraction. We initially prove that $\Phi_{\alpha}$ is well defined. By the expression (3.4), we can claim that $S_{\alpha}(\cdot)[\varphi(0)-g(0,\varphi)]$ belongs to $PSAP_{\omega,r}(X)$. Moreover, for $\forall~u\in \mathcal{B}$ and $\forall~t\geq0$, we have
\begin{eqnarray*}
\|\Phi_{\alpha}u(t)\|&\leq&\|S_{\alpha}(t)\|[\|\varphi(0)\|+\|g(0,\varphi)\|]+\bigg[\|g(t, u_{t})-g(t,0)\|+\|g(\cdot,0)\|_{C_b([0,\infty),X)}\bigg]\\
&&+\int_{0}^{t}\frac{CM}{1+|\mu|(t-s)^{\alpha}}\bigg[\|f(s, u_{s})-f(s,0)\|+\|f(\cdot,0)\|_{C_b([0,\infty),X)}\bigg]\\
&\leq&CM[\|\varphi(0)\|+\|g(0,\varphi)\|]+\bigg[\|L_g\|_{C_b([0,\infty),\mathbb{R}^{+})}\|u\|_{C_b([-r,\infty),X)}+\|g(\cdot,0)\|_{C_b([0,\infty),X)}\bigg]\\
&&+\sup\limits_{t\geq0}W_{f}(t)\|u\|_{C_b([-r,\infty),X)}+\displaystyle\frac{CM|\mu|^{\frac{-1}{\alpha}}\pi}{\alpha\sin(\frac{\pi}{\alpha})}\|f(\cdot,0)\|_{C_b([0,\infty),X)}.
\end{eqnarray*}
Combing with Lemmas 2.4-2.5 and the arguments above, it follows that $\Phi_{\alpha}u$ is well defined on $\mathcal{B}$.

Furthermore, for all $u,v\in \mathcal{B}$ and $t\geq0$, we have
\begin{eqnarray*}
&&\|(\Phi_{\alpha}u)(t)-(\Phi_{\alpha}v)(t)\|\\
&&~\leq L_{g}(t)\|u_t-v_t\|_{\mathcal{C}}+CM\int_{0}^{t}\frac{L_f(s)}{1+|\mu|(t-s)^{\alpha}}\|u_s-v_s\|_{\mathcal{C}}ds\\
&&~\leq\bigg[\|L_{g}\|_{C_b([0,\infty),\mathbb{R}^{+})}+\sup\limits_{t\geq0}W_{f}(t)\bigg]\|u-v\|_{C_b([-r,\infty),X)}.
\end{eqnarray*}
The assumption (3.5) enables us to claim our assertion. This completes the proof of Theorem 3.3.$\quad\Box$\\[3mm]
{\bf Theorem 3.4.} \emph{Assume that hypotheses} (H$_1$) \emph{and} (H$_5$) \emph{hold, and
$$\|L_g\|_{C_b([0,\infty),\mathbb{R}^{+})}+CM\bigg(1+\displaystyle\frac{|\mu|^{\frac{-1}{\alpha}}\pi}{\alpha\sin(\frac{\pi}{\alpha})}\bigg)\|L_f\|_{S^p}<1,\eqno(3.6)$$
then system $(1.2)$ has a unique mild solution $u\in PSAP_{\omega,r}(X)$.}\\[3mm]
\emph{ Proof.} The operator $\Phi_{\alpha}$ and the set $\mathcal{B}$ are defined as in Theorem 3.1. From the assumption (H$_5$), according to Theorem 3.2 in [37], we know that $f(s,u_s)\in PSAP_{\omega,r}(X)$. Proceeding as in the proof of Theorem 3.3, it is easy to claim that $\Phi_{\alpha}u\in \mathcal{B}$ whenever $u\in \mathcal{B}$.\\
Furthermore, for all $u, v \in \mathcal{B}$ and each $t\geq0$, one has
\begin{eqnarray*}
&&\|(\Phi_{\alpha}u)(t)-(\Phi_{\alpha}v)(t)\|\\
&&~\leq L_{g}(t)\|u_t-v_t\|_{\mathcal{C}}+\int_{0}^{t}S_{\alpha}(t-s)\|f(s,u_s)-f(s,v_s)\|ds\\
&&~\leq L_{g}(t)\|u_t-v_t\|_{\mathcal{C}}+CM\int_{0}^{t}\frac{L_f(s)}{1+|\mu|(t-s)^{\alpha}}\|u_s-v_s\|_{\mathcal{C}}ds.
\end{eqnarray*}
 If $t=m \in \mathbb{N}$, in this case
\begin{eqnarray*}
&&\|(\Phi_{\alpha}u)(t)-(\Phi_{\alpha}v)(t)\|\\
&&~\leq L_{g}(t)\|u_t-v_t\|_{\mathcal{C}}+CM\int_{0}^{m}\frac{L_f(s)}{1+|\mu|(m-s)^{\alpha}}\|u_s-v_s\|_{\mathcal{C}}ds\\
&&~\leq L_{g}(t)\|u_t-v_t\|_{\mathcal{C}}+CM\sum\limits_{k=0}^{m-1}\int_{k}^{k+1}\frac{L_f(s)}{1+|\mu|(m-s)^{\alpha}}\|u_s-v_s\|_{\mathcal{C}}ds\\
&&~\leq L_{g}(t)\|u_t-v_t\|_{\mathcal{C}}+CM\sum\limits_{k=0}^{m-1}\frac{1}{1+|\mu|(m-k-1)^{\alpha}}\sup\limits_{k\in[0,m-1]}\int_{k}^{k+1}L_f(s)ds\|u_s-v_s\|_{\mathcal{C}}\\
&&~\leq L_{g}(t)\|u_t-v_t\|_{\mathcal{C}}+CM\bigg[1+(\int_{0}^{1}+\int_{1}^{2}+\cdot\cdot\cdot+\int_{m-2}^{m-1})\frac{1}{1+|\mu|\tau^{\alpha}}d\tau\bigg]\\
&&~~~\times\bigg(\sup\limits_{k\in[0,m-1]}\int_k^{k+1}\|L_f(s)\|^pds\bigg)^{\frac{1}{p}}\|u_s-v_s\|_{\mathcal{C}}\\
&&~\leq L_g(t)\|u_t-v_t\|_{\mathcal{C}}+CM\bigg[1+\int_{0}^{\infty}\frac{1}{1+|\mu|\tau^{\alpha}}d\tau\bigg]\|L_f\|_{S^p}\|u_s-v_s\|_{\mathcal{C}}\\
&&~\leq\bigg[\| L_g\|_{C_b([0,\infty),\mathbb{R}^{+})}+CM+\frac{CM|\mu|^{\frac{-1}{\alpha}}\pi}{\alpha\sin(\frac{\pi}{\alpha})}\|L_f\|_{S^p}\bigg]\|u-v\|_{C_b([-r,\infty),X)}.
\end{eqnarray*}
We now deal with more general case. If $t = m - h$, where $0 < h < 1,~m\in \mathbb{N}$, we have
\begin{eqnarray*}
&&\int_{0}^{t}S_{\alpha}(t-s)L_f(s)d\leq\int_{0}^{m-h}\frac{CM}{1+|\mu|(m-h-s)^{\alpha}}L_f(s)ds\\
&&~\leq\int_{h}^{m}\frac{CM}{1+|\mu|(m-s)^{\alpha}}L_f(s-h)ds
\leq \int_{h}^{m}\frac{CM}{1+|\mu|(m-s)^{\alpha}}\bar{L}_f(s)ds\\
&&~\leq \int_{0}^{m}\frac{CM}{1+|\mu|(m-s)^{\alpha}}\bar{L}_f(s)ds
\leq (\int_{0}^{h}+\int_{h}^{m})\frac{CM}{1+|\mu|(m-s)^{\alpha}}\bar{L}_f(s)ds\\
&&~\leq CM\bigg (1+\frac{|\mu|^{\frac{-1}{\alpha}}\pi}{\alpha\sin(\frac{\pi}{\alpha})}\bigg)\|\bar{L}_f\|_{S^p}.
\end{eqnarray*}
where $\bar{L}_f$ is defined by
$$
 \bar{L}_f(s):=\left\{\begin{array}{ll}
 0,~~~~~~~~~0\leq s< h,\\
  L_f(s-h),~~~s\geq h.
  \end{array}
  \right.$$

Thus we have $\|\bar{L}_f\|_{S^p}=\|L_f\|_{S^p}$. Moreover, we can easily get
$$\int_{0}^{t}\frac{1}{1+|\mu|(t-s)^{\alpha}}L_f(s)ds\leq\bigg(1+\frac{|\mu|^{\frac{-1}{\alpha}}\pi}{\alpha\sin(\frac{\pi}{\alpha})}\bigg)\|\bar{L}_f\|_{S^p}.$$ From the preceding estimates, we have
$$\|\Phi_{\alpha}u-\Phi_{\alpha}v\|_{C_b([-r,\infty),X)}\leq\bigg[\|L_g\|_{C_b([0,\infty),\mathbb{R}^{+})}+CM\bigg(1+\frac{|\mu|^{\frac{-1}{\alpha}}\pi}{\alpha\sin(\frac{\pi}{\alpha})}\bigg)\|L_f\|_{S^p}\bigg]\|u-v\|_{C_b([-r,\infty),X)}.$$
By (3.6) and the Banach contraction mapping principle, we know that $\Phi_{\alpha}$ has a unique fixed point $u\in PSAP_{\omega,r}(X)$ in $\mathcal{B}$, which is the mild solution of system (1.2). The proof is complete.$\quad\Box$

%%%%%%%%%%%%%
% SECTION 5 %
%%%%%%%%%%%%%

\section{Applications}
\setcounter{equation}{0}\setcounter{section}{4}\indent
In what follows, we use the previous theory to verify the existence and uniqueness of pseudo $S$-asymptotically
$\omega$-periodic mild solutions to a fractional partial integro-differential neutral equation. We are concerned with the following system
$$
\left \{\begin{array}{lll}
\displaystyle\frac{\partial}{\partial t}\bigg[v(t,\vartheta)-h(t)\int_{-r}^{0}k(s)v(t+s,\vartheta) ds\bigg]=J_t^{\alpha-1}(\frac{\partial^2}{\partial\vartheta}-c)\bigg[v(t,\vartheta)-h(t)\int_{-r}^{0}k(s)v(t+s,\vartheta)ds\bigg]\\
~~~~~~~~~~~~~~~~~~~~~~~~~~~~~~~~~~~~~~~~~~~~~~~+j(t)\displaystyle\int_{-r}^{0}m(s)v(t+s,\vartheta)ds,~~t\geq0,~c>0,~\vartheta\in[0,\pi],\\
 v(\theta,\vartheta)=\varphi(\theta,\vartheta),~~~~~\theta\in[-r,0],
\end{array}
\right.\eqno(4.1)
$$
where $1<\alpha<2,~h(t),j(t)\in PSAP_{\omega,r}([0,+\infty),\mathbb{R})$, $k,m\in C_b([-r,0],\mathbb{R})$ and satisfy some  particular conditions specified later. Let $X=L^2([0,\pi],\mathbb{R},\|\cdot\|_{L^2})$ and define the operator $A$ on $X$ by
$Av=(\displaystyle\frac{\partial^2}{\partial\vartheta}-c)v$. $D(A)=\{v\in X=L^{2}([0,\pi],\mathbb{R}):v^{''}\in L^2[0,\pi],~v(0)=v(\pi)=0\}$. It is well known that $A$ is sectorial of type $\omega = -c < 0$. Let $u(t)(\vartheta)=v(t,\vartheta)$ for $t\in[-r,\infty),$ $g(t,u_t)(\vartheta)=h(t)\displaystyle\int_{-r}^{0}k(s)v(t+s,\vartheta) ds,$ $f(t,u_t)(\vartheta)=j(t)\displaystyle\int_{-r}^{0}m(s)v(t+s,\vartheta) ds,~t\in [0,\infty),~\vartheta\in [0,\pi],$ and $J_t^{\alpha-1}u(t)=\displaystyle\int_{0}^{t}\frac{(t-s)^{\alpha-2}}{\Gamma(\alpha-1)}u(s)ds.$ Then (4.1) can be rewritten as an
abstract system of the form (1.2).

For reader$^{,}$s convenience, we only give a brief outline of the proof. We have the following estimates
$$
\begin{array}{rl}
\|g(t,\varphi)\|\leq r^{\frac{1}{2}}\|h\|\bigg(\displaystyle\int_{-r}^{0}|k(s)|^2ds\bigg)^{\frac{1}{2}}\|\varphi\|_{\mathcal{C}},
\end{array}\eqno(4.2)
$$
$$
\begin{array}{rl}
\|g(t,\varphi)-g(t,\psi)\|\leq r^{\frac{1}{2}}\|h\|\bigg(\displaystyle\int_{-r}^{0}|k(s)|^2ds\bigg)^{\frac{1}{2}}\|\varphi-\psi\|_{\mathcal{C}},\\
\end{array}\eqno(4.3)
$$
$$
\begin{array}{rl}
&\displaystyle\frac{1}{T}\int_{p}^{T}\sup\limits_{\tau\in[t-r,t]}\|g(\tau+\omega,u_{\tau})-g(\tau,u_{\tau})\|dt\\
&\leq\displaystyle\frac{1}{T}\int_{p}^{T}\sup\limits_{\tau\in[t-r,t]}\|h(\tau+\omega)-h(\tau)\|dt\times\|r^{\frac{1}{2}}\|(\int_{-r}^{0}k^{2}(s)ds)^{\frac{1}{2}}\|u_{\tau}\|_{\mathcal{C}}.
\end{array}\eqno(4.4)
$$
 We note that estimate (4.2) implies that $g$ is bounded. Since $h$ belongs to $PSAP_{\omega,r}(X)$, then estimates (4.3) and (4.4) imply that $g\in PSAP_{\omega,r}(\mathcal{C},X)$ with $L_g=r^{\frac{1}{2}}\|h\|\bigg(\displaystyle\int_{-r}^{0}|k(s)|^2ds\bigg)^{\frac{1}{2}}$. Similarly, we can verify that $f$ is bounded and $f\in PSAP_{\omega,r}(\mathcal{C},X)$ with $L_f=r^{\frac{1}{2}}\|j\|\bigg(\displaystyle\int_{-r}^{0}|m(s)|^2ds\bigg)^{\frac{1}{2}}$. Thus we can deduce the following result as an immediate consequence
of Corollary 3.2.\\[3mm]
{\bf Proposition 4.1.} \emph{Under previous assumptions, if we assume $L_f$ and $L_g$ are constants small enough such that $L_g+L_f\displaystyle\frac{CM|\mu|^{\frac{-1}{\alpha}}\pi}{\alpha\sin(\frac{\pi}{\alpha})}<1,$ where $L_f$ and $L_g$ are introduced above. Then there exists a unique  mild
solution $v\in PSAP_{\omega,r}(X)$ of system (4.1).}

\Acknowledgements{This work was supported by National Natural Science Foundation of China (Nos. 11271379 and 11671406). We would like to thank the reviewers, editor and associate
editor-in-chief for their valuable suggestions and comments.}



\begin{thebibliography}{99}
\bahao\baselineskip 11.5pt

\bibitem{Abbas2011}
Abbas S.
 Pseudo almost automorphic solutions of some nonlinear integro-differential equations. Comput Math Appl, 2011, 62: 2259-2272

 \bibitem{Agarwal2010}
Agarwal R P, Andrade B de, Cuevas C.
Weighted pseudo-almost periodic solutions of a class of semilinear
fractional differential equations. Nonlinear Anal Real World Appl, 2010, 11: 3532--3554


\bibitem{Agarwal2011}
Agarwal R P, Cuevas C, Soto H, El-Gebeily M.
Asymptotic periodicity for some evolution equations in Banach spaces. Nonlinear Anal, 2011, 74:  1769--1798 %


\bibitem{Alvarez2015}
Alvarez E, Lizama C, Ponce R.
Weighted pseudo antiperiodic solutions for fractional
integro-differential equations in Banach spaces. Appl Math Comput, 2015, 259: 164--172 %


\bibitem{Andrade2010}
 Andrade B de, Cuevas C.
$S$-asymptotically $\omega$-periodic and asymptotically $\omega$-periodic solutions to semilinear Cauchy problems with nondense domain. Nonlinear Anal, 2010, 72: 3190--3208 %


\bibitem{Andrade2015}
Andrade F, Cuevas C, Silva C, Soto H.
Asymptotic periodicity for hyperbolic evolution equations and applications. Appl Math Comput, 2015, 269: 169--195 %


\bibitem{Bochner1962}
Bochner S.
 A new approach to almost periodicity. Proc Natl Acad Sci USA, 1962, 48: 2039--2043 %


\bibitem{Caicedo2012}
Caicedo A, Cuevas C, Mophou G M, \'{N}Gu\'{e}r\'{e}kata G M.
Asymptotic behavior of solutions of some semilinear functional differential and integrodifferential
equations with infinite delay in Banach spaces. J Franklin Inst, 2012, 349: 1--24 %


\bibitem{Agarwal2010}
Agarwal R P, Andrade B de, Cuevas C.
On type of periodicity and ergodicity to a class of fractional order differential equations. Adv Difference Equ,
2010, Article ID 179750, 25 pp.


\bibitem{Chen2011}
 Chen X X, Hu X Y.
   Weighted pseudo almost periodic solutions of neutral functional
differential equations. Nonlinear Anal Real World Appl, 2011, 12: 601--610 %


\bibitem{Corduneanu1989}
Corduneanu C.
 Almost Periodic Functions, New York: Wiley, 1968; Reprinted, New York: Chelsea, 1989 %


\bibitem{Cuesta2007}
Cuesta E.
Asymptotic behavior of the solutions of fractional integro-differential equations and some time discretizations. Discrete Contin Dyn
Syst, (Supplement) 2007, 277--285 %


\bibitem{Cueva2009}
Cuevas C, Hern\'{a}ndez E.
 Pseudo-almost periodic solutions for abstract partial functional differential equations. Appl Math Lett, 2009, 22: 534--538 %


\bibitem{Cuevas2014}
Cuevas C, Henr\'{\i}quez H R, Soto H.
Asymptotically periodic solutions of fractional differential equations. Appl Math Comput, 2014, 236: 524--545 %


\bibitem{Cuevas2012}
 Cuevas C, de Souza J C.
 $S$-asymptotically $\omega$-periodic solutions of semilinear fractional integro-differential equations. Appl Math Lett, 2009, 22: 865--870 %


\bibitem{Cuevas2010}
Cuevas C, de Souza J C.
 Existence of $S$-asymptotically $\omega$-periodic solutions for fractional order functional integro-differential equations with infinite delay. Nonlinear Anal, 2010, 72: 1683--1689 %


\bibitem{Cuevas2011}
Cuevas C, Lizama C.
$S$-asymptotically $\omega-$periodic solutions for semilinear Volterra equations. Math Methods Appl Sci, 2010, 33: 1628--1636 %


\bibitem{Cuevas2013}
Cuevas C, Lizama C.
Existence of $S$-asymptotically $\omega$-periodic solutions for two-times fractional order differential equations. Southeast Asian Bull Math, 2013, 37: 683--690 %


\bibitem{Diagana2013}
Diagana T.
  Almost Automorphic Type and Almost Periodic Type Functions in Abstract Spaces. New York: Springer, 2013 %


\bibitem{dos Santos2010}
dos Santos J P C, Cuevas C.
Asymptotically almost automorphic solutions of abstract fractional integro-differential neutral equations. Appl Math Lett,
2010, 23: 960--965 %

\bibitem{dos Santos2012}
dos Santos J P C, Rabelo M, Henrique M.
Asymptotic almost automorphic and $S$-asymptotically $\omega$-periodic solutions to partial differential equations. Adv
Differ Equ Control Process, 2012, 9: 45--61 %


\bibitem{Rolnik2016}
Henr\'{\i}quez H R, Pierri M, Rolnik V.
  Pseudo $S$-asymptotically periodic solutions of second-order abstract Cauchy problems. Appl Math Comput, 2016, 274: 590--603 %


 \bibitem{Henriquez2008}
Henr\'{\i}quez H R, Pierri M, T\'{a}boas P.
On $S$-asymptotically $\omega$-periodic functions on Banach spaces and applications. J Math Anal Appl, 2008, 343: 1119--1130 %


\bibitem{Henriquez2008}
 Henr\'{\i}quez H R, Pierri M, T\'{a}boas P.
Existence of $S$-asymptotically $\omega$-periodic solutions for abstract neutral equations. Bull Aust Math Soc, 2008, 78: 365--382 %


\bibitem{Ji2012}
Ji D S, Zhang C Y.
Translation invariance of weighted pseudo almost periodic functions and
related problems. J Math Anal Appl, 2012, 391: 350--362 %

\bibitem{Kilbas2006}
 Kilbas A, Srivastava H, Trujillo J.
  Theory and Applications of Fractional Differential Equations, North-Holland Mathematics Studies 204, Amsterdam: Elsevier Science B.V, 2006 %

\bibitem{Li2014}
 Li Z, Liu K, Luo J W.
  On almost periodic mild solutions for neutral stochastic
evolution equations with infinite delay. Nonlinear Anal, 2014, 110: 182--190 %


 \bibitem{Mophou2011}
Mophou G M.
 Weighted pseudo almost automorphic mild solutions to semilinear fractional differential equations. Appl Math Comput, 2011, 217:
7579--7587 %


\bibitem{NGuerekata2005}
N'Gu\'{e}r\'{e}kata G M.
 Topics in Almost Automorphy, New York: Springer, 2005 %


\bibitem{Nicota2009}
Nicota S, Pierri M.
A note on $S$-asymptotically periodic functions. Nonlinear Anal Real World Appl, 2009, 10: 2937--2938 %


\bibitem{Pierri2012}
Pierri M,
On $S$-asymptotically $\omega$-periodic functions on Banach spaces and applications. Nonlinear Anal, 2012, 75: 651--661 %


\bibitem{Pierri2013}
Pierri M, Rolnik V.
On pseudo $S$-asymptotically $\omega$-periodic functions. Bull Aust Math Soc, 2013, 87: 238--254 %

\bibitem{Revathi2014}
Revathi P, Sakthivel R, Ren Y, Anthoni S M.
Existence of almost automorphic mild solutions to
non-autonomous neutral stochastic differential equations. Appl Math Comput, 2014, 230: 639--649 %

\bibitem{Sakthivel2012}
Sakthivel R, Revathi P, Anthoni S M.
Existence of pseudo almost automorphic mild solutions to stochastic
fractional differential equations. Nonlinear Anal, 2012, 75: 3339--3347 %


\bibitem{Shu2015}
Shu X B,  Xu F, Shi Y J.
$S$-asymptotically $\omega$-positive periodic solutions for a class of neutral fractional differential equations. Appl Math Comput, 2015, 270: 768--776 %


\bibitem{Xia2014}
 Xia Z N,
Weighted pseudo periodic solutions of neutral functional differential equations, Electron J Differ Eq, 2014, 191: 1--17 %


\bibitem{Xia2015}
Xia Z N,
Pseudo asymptotically periodic solutions of two-term time fractional differential equations with delay. Kodai Math J, 2015, 38: 310--332 %


\bibitem{Zhang2012}
Zhang R, Chang Y K, \'{N}Gu\'{e}r\'{e}kata G M.
 New composition theorems of Stepanov-like almost automorphic functions and applications to
nonautonomous evolution equations. Nonlinear Anal Real World Appl, 2012, 13: 2866--2879 %


\bibitem{Zhang2015}
Zheng Z M, Ding H S.
On completeness of the space of weighted pseudo almost automorphic functions. J Funct Anal, 2015, 268: 3211--3218 %


\bibitem{Zhou2014}
Zhou Y.
 Basic Theory of Fractional Differential Equations. New Jersey: World Scientific, 2014 %


\bibitem{Zhou2010}
Zhou Y, Jiao F.
 Existence of mild solutions for fractional neutral evolution equations. Comput Math Appl, 2010, 59: 1063--1077

\end{thebibliography}
\end{document}